\documentclass[10pt]{amsart}
\usepackage{amsmath,amssymb,latexsym,amsthm}
\font\goth=eusm10
\newcommand\E{\mathcal E}
\newcommand\oE{\mathcal {E'}}
\newcommand\Ii{\hbox{\goth I}}

\newcommand\Oc{\hbox{\goth O}}

\newcommand\CC{\mathbb C}
\newcommand\Hh{H}

\newcommand\Pj{\mathbb{P}}

\newcommand\Pt{\mathbb{P}^3}
\newcommand\Pq{\mathbb{P}^4}
\newcommand\Pc{\mathbb{P}^5}
\newcommand\ZZ{\mathbb{Z}}
\newtheorem{theorem}{Theorem}[section]
\newtheorem{proposition}[theorem]{Proposition}
\newtheorem{definition}[theorem]{Definition}
\newtheorem{corollary}[theorem]{Corollary}

\newtheorem{remark}[theorem]{Remark}

\newtheorem{case}{Case}[section]
\theoremstyle{plain}
\theoremstyle{definition}
\theoremstyle{remark}
\numberwithin{equation}{section}

\begin{document}
\title
{ACM bundles on general hypersurfaces in $\Pj^5$ of low degree}
\author{L.Chiantini, C.Madonna}
\address{L.Chiantini, Dipartimento di Scienze Matematiche e Informatiche,
Universit\`a di Siena, Pian dei Mantellini 44, 53100 SIENA  (Italy)}
\email{chiantini@unisi.it}
\address{C.K.Madonna, Dipartimento di Matematica, Universit\`a degli Studi di 
Roma "La Sapienza",  P.le A.Moro 1, 00100 ROMA (Italy)}
\email{madonna@mat.uniroma1.it}

\begin{abstract}
In this paper we show that on a general hypersurface of degree $r=3,4,5,6$ in $\Pj^5$
a rank $2$ vector bundle $\E$ splits if and only if $h^1\E(n)=h^2 \E(n)=0$ for all 
$n \in \ZZ$. Similar results for $r=1,2$ were obtained in \cite{OSS}, \cite{Ott} 
and \cite{arrondo}.
\end{abstract}
\maketitle

\section{Introduction} \label{S:intro}

The construction of rank $2$ bundles on smooth varieties $X$ of dimension $n> 3$ is 
strictly related with the structure of subvarieties of codimension $2$. When $X$ is a 
projective space, then there are few examples of these subvarieties which are smooth. 
The famous Hartshorne's conjecture suggests that all smooth subvarieties of codimension 
$2$ in $\Pj^7$ are complete intersection. Rephrased in the language of vector bundles, 
this means that all rank $2$ bundles on $\Pj^7$ decompose in a sum of two line bundles.\par
Also in $\Pj^5$, $\Pj^6$, we do not have examples of indecomposable rank $2$ bundles. 
In $\Pj^4$, only the Horrocks-Mumford's indecomposable bundle is known. This bundle has 
some non--zero cohomology group, since it is well known that a rank $2$ bundle $\E$ on 
$\Pj^r$ ($r\geq 3$) splits if and only if it is ``arithmeticall Cohen--Macaulay'' (ACM 
for short), i.e. $h^i(\E(n))=0$ for all $n, i=1,\dots,r-1$.\par
ACM property does not imply a decomposition when we replace the projective space with 
other smooth threefolds. There are examples of indecomposable ACM bundles of rank $2$ 
on a general hypersurface of degree $r=2,3,4,5$ in $\Pq$. On the other hand we proved 
in \cite{cm2} that all ACM rank $2$ bundles on a {\it general} sextic in $\Pq$ splits.\par
\smallskip
In this paper we examine the similar problem for general hyperpsurfaces $X$ in $\Pc$, 
in some sense the easiest examples of smooth $4$-folds different from $\Pq$.\par   
It is well known that a general quadric hypersurface $X$ in $\Pc$ contains families of 
planes. Since any plane $S$ has a canonical class which is a twist of the restriction of 
the canonical class of the quadric (in other words: a plane is ''subcanonical'' in $X$), 
then $S$ corresponds via the Serre's construction to a rank $2$ bundle $\E$ on $X$ which 
is indecomposable (for $S$ is not complete intersection of $X$ and some other 
hypersurface) and ACM (for $S$ is arithmetically Cohen--Macaulay).\par
On the other hand, since any indecomposable ACM rank $2$ bundle on a general sextic 
hypersurface in $\Pc$ would restrict to an indecomposable rank $2$ ACM bundle on a 
general hyperplane section of $X$, which is a general sextic hypersurface of $\Pq$, 
then by the main result of \cite{cm2} we know that such bundles cannot exist (see 
proposition \ref{prop:sextic} below).\par
Thus we are led to consider {\it general} hypersurfaces $X_r\subset\Pc$ of low degree 
$r$ and study ACM rank $2$ bundles on $X_r$. Our main result shows that none of such 
vector bundles lives on $X_r$  for $2<r<7$:

\begin{theorem} \label{thm:main}
Let $\E$ be a rank $2$ vector  bundle on a general hypersurface $X_r \subset \Pj^5$
of degree $r=3,4,5,6$. Then $\E$ splits if and only if  
$$h^i(\E(n))=0 \qquad \forall  n \in \ZZ \qquad i=1,2.$$
\end{theorem}

Notice that one finds indecomposable ACM rank $2$ bundles on general hypersurfaces of 
degree $3,4,5$ in $\Pq$. So we prove in fact that they do not lift from  a general 
hyperplane section of $X$ to $X$ itself.\par
The proof is achieved using the tools of \cite{cm2}, since we have a classification of 
possible ACM indecomposable rank $2$ bundles on a general hypersurface of low degree in 
$\Pq$. It has been obtained by Arrondo and Costa in degree $3$ (see \cite{AC}), by the 
second author in degree $4$ (see \cite{m2}) and by both authors in degree $5$ 
(\cite{cm1}). This implies a numerical characterization of the possible Chern classes of 
indecomposable ACM bundles of rank $2$ on $X$ (see also \cite{m0}), and we conclude with 
a case by case analysis.\par
In the language of codimension $2$ subvarieties, we get the following characterization 
of complete intersections, which is the analogue of the classical Gherardelli's criterion 
for curves in $\Pt$:

\begin{corollary}
Let $S$ be a surface contained in a general  hypersurface $X_r \subset \Pj^5$ of degree 
$r= 3,4,5,6$. Then $S$ is complete intersection in $X_r$ if and only if it is subcanonical 
(i.e. its canonical class $\omega_X$ is $\Oc_S(e)$, for some $e\in\ZZ$) and 
$h^i \Ii_{S/X_r}(n)=0$ for all $n \in \ZZ$ and $i=1,2$, where $\Ii_{S/X_r}$ is the ideal 
sheaf defining $S$ in  $X_r$.
\end{corollary}

Let us finish with some remarks.\par\smallskip

The non--existence of indecomposable ACM rank $2$ bundles on hypersurfaces of degree 
$r\geq 7$ in $\Pq$ has not been settled yet simply because the tecnicalities introduced 
in \cite{cm2} become odd as the degree $r$ grows. Indeed first of all the number of
Chern classes which are not excluded using the main result of \cite{m0} grows as a linear
function of $r$. Furthermore, as $r$ grows, for any value of $c_1$ one has to exclude
an increasing number of second Chern classes. This is easy when $c_2$ is big, but becomes
hard for low $c_2$ (compare the proof of case 5.11 in \cite{cm2}),
as we have to exclude the existence of some curves on $X$, which could be reducible
or even non--reduced. We did not find a general argument for this step: only a careful 
ad hoc examination led us to conclude the case of sextic threefolds in $\Pq$ .\par  

On the other hand, there are strong evidences that ACM rank $2$ indecomposable 
bundles cannot exist on general hypersurfaces of degree $7$ or more. We were unable 
to prove this statement in $\Pq$. Could it be easier to find a  direct proof for 
hypersurfaces in $\Pc$? 
\par\smallskip

In any event, the main theorem implies easily:

\begin{corollary} On a general hypersurface $X$ of degree $3,4,5,6$ in $\Pj^n$, 
$n\geq 5$, a rank $2$ vector bundle splits if and only if it is arithmetically 
Cohen--Macaulay.
\end{corollary}

Finally observe that Evans and Griffith proved in \cite{EG} that a rank $2$ bundle 
$\E$ on $\Pq$ splits if and only $h^1(\E(n))=0$ for all $n$. This condition is 
considerably weaker than ACM. We wonder if a similar result could work on a general 
hypersurface of low degree in $\Pc$.

\section{Notations and generalities} \label{S:gener}

We work over the field of complex numbers $\CC$.
Let $X_r \subset \Pc$ be a smooth $4$--dimensional hypersurface of degree $r \geq 1$. 
The letter $H$ will denote the class of a hyperplane section of $X_r$. We have $Pic(X_r) 
\cong \ZZ[H]$, and $H^4=r$. Recall that the canonical class of $X_r$ is
$\omega_{X_r}=(r-6)H$. Given a vector bundle $\E$ on $X_r$ we introduce the non 
negative integer
\begin{equation} \label{eq:norm}
b(\E)=b=\max \{ n \mid h^0(\E(-n)) \ne 0 \}.
\end{equation}

\begin{definition}
We say that the vector bundle $\E$ is normalized if $b(\E)=0$.
\end{definition}

Notice that changing $\E$ with $\E(-b)$, we may always assume that $\E$ 
is normalized. From now on we will assume this.

We denote by $c_1=c_1(\E)$ the first Chern class of $\E$ identified with an integer via 
the isomorphism $Pic(X_r) \cong \ZZ[H]$.

When $\E$ has rank $2$, the number 
\begin{equation} \label{eq:2b-c1}
2b-c_1 = 2b(\E)-c_1(\E)
\end{equation} 
is invariant by twisting and measures the ``level of stability of $\E$''. 
Indeed $\E$ is stable (semistable) if an only if $0>2b-c_1$ ($0 \geq 2b-c_1$).

We say that $\E$ is ``arithmetically Cohen--Macaulay (ACM)'' when for 
all $n \in \ZZ$ we have $h^1(\E(n))=h^2(\E(n))=0$. Clearly this implies, by duality, 
$h^3(\E(n))=0$ for all $n \in \ZZ$. 

Take a global section $s$ of $\E$ whose zero-locus $S$ has codimension $2$. We have the 
following exact sequence (see e.g. \cite {v}):
\begin{equation} \label{eq:kos}
0 \to \Oc_{X_r} \to \E \to \Ii_{S/X_r}(c_1(\E)) \to 0
\end{equation}
which relates the cohomology of $\E$ with the geometric properties of $S \subset X_r$ 
encoded by the cohomology groups of the ideal sheaf $\Ii_{S/X_r}$ of $S$.

In particular $S$ is subcanonical, i.e. $K_S \cong \Oc_S(c_1(\E)+r-6)$, moreover $
c_2(\E)=\deg S$. Also we have a formula for the sectional genus $g$ of the surface $S$:
\begin{equation} \label{eq:sectional}
2g-2=c_2+K_S \cdot H \cdot S=c_2+(c_1+r-6)H \cdot H \cdot S=
c_2(c_1+r-5)
\end{equation}

Conversely, starting with a locally complete intersection and subcanonical surface $S$ 
contained in $X_r$ one can reconstruct a rank 2 vector bundle having a global section 
vanishing exactly on $S$. In these cases we will say that the vector bundle ``$\E$ is 
associated with $S$''.

We notice that when $\E$ is normalized, then every global section of $\E$ has zero-locus  
of codimension $2$.

If $Y_r$ is a general hyperplane section of $X_r$ and $\E$ is a rank two vector bundle 
on $X_r$, we denote by $\oE$ the restriction of $\E$ to $Y_r$. We know that 
$Pic(Y_r)\cong Z[h]$, where $h$ is the hyperplane class of $Y_r$. Under the isomorphism 
$Pic(X_r) \cong Pic(Y_r)$ we have $c_1(\E)=c_1(\oE)$ and $c_2(\E)=c_2(\oE)$.

We recall here the main results of \cite{m0} and \cite{cm2}, which we are going to
use several times in the sequel:

\begin{theorem} \label{m0} (see \cite{m0})
Let $Y_r$ be a smooth hypersurface of degree $r$ in $\Pq$. If $\E$ is an ACM and 
normalized rank $2$ vector bundle on $Y_r$, then $\E$ splits unless $r> c_1 > 2-r$.
\end{theorem}

\begin{theorem} \label{cm2} (see \cite{cm2})
Let $Y$ be a general hypersurface of degree $6$ in $\Pq$. Then a rank $2$ vector bundle
$\E$ on $Y$ splits in a sum of line bundles if and only if $\E$ is ACM.
\end{theorem}

\section{Some preliminary general results}

\begin{remark}\label{rm:restrict} 
\rm Consider the exact sequence which links $\E$ with its restriction $\oE$ to a general 
hyperplane section $Y_r$ of $X_r$:
\begin{equation} \label{seq:restr}
0 \to \E(-1) \to \E \to \oE \to 0
\end{equation}

Then $b(\oE)\geq b(\E)$ and equality holds when $h^1(\E(-b(\E)-2))=0$, which is true 
when $\E$ is ACM.

Notice that, by the sequence, if $\E$ is ACM on $X_r$ then $\oE$ is also ACM on $Y_r$.
It is clear that $\oE$ splits when $\E$ splits. 
Conversely assume that $\E$ is ACM and $\oE$ splits. Take a global section $s'\in 
H^0(\oE(a))$ with empty zero-locus. The surjection $H^0(\E(a)) \to H^0(\oE(a)) \to 0$ 
derived from sequence  (\ref{seq:restr}) shows that $s'$ lifts to a global section 
$s\in H^0(\E(a))$, whose zero-locus must be empty, since otherwise it had at most 
codimension $2$, a contradiction for it does not intersect a hyperplane.\end{remark}

It follows that we may apply the main result of \cite{m0}, getting:

\begin{proposition} \label{prop:m0}
If $\E$ is an ACM and normalized rank $2$ vector bundle on $X_r$,
then $\E$ splits unless $r> c_1 > 2-r$.
\end{proposition}

Some well known facts about the non--existence of surfaces of low degree on general 
hypersurfaces of $\Pc$ together with a numerical analysis leads us to the following 
refinement of the previous result:

\begin{proposition} \label{prop:refin}
Let $\E$ be a normalized ACM rank 2 vector bundle on a general 
hypersurface $X_r\subset \Pc$ of degree $r \geq 3$. Then $\E$ splits unless 
$3-r < c_1 < r$. 
\end{proposition}
\begin{proof}
We need to exclude the case $c_1=3-r$.\par\noindent
Consider a global section $s$ and its zero-locus $S$. The exact sequence 
(\ref{eq:kos}) here reads  
$$
0 \to \Oc_{X_r} \to \E \to \Ii_{S/X_r}(3-r) \to 0
$$ 
and implies $h^0(\E(r-3))= h^0(\Oc_{X_r}(r-3))$. By Serre 
duality $h^4(\E(r-3))=h^0(\E(r-6))=h^0(\Oc_{X_r}(r-6))$. Moreover $h^1(\E(r-3))=
h^2(\E(r-3))=0$ for $\E$ is ACM. Thus $\chi(\E(r-3))=h^0(\Oc_{X_r}(r-3))+
h^0(\Oc_{X_r}(r-6))$. By Riemann-Roch one is thus able to compute the second Chern class 
of $\E(r-3)$, hence also the second Chern class $c_2$ of $\E$. It turns out $c_2=1$. 
So $S$ is a plane. Since a general $X_r$ of degree $r \geq 3$ contains no planes (see 
e.g. \cite{debarre}), then $X_r$ has no indecomposable and normalized rank 2 ACM bundles 
with $c_1=-r+3$.
\end{proof}

Next we use the link between ACM bundles with $c_1=r-1$ 
and pfaffian hypersurfaces.

\begin{definition} A hypersurface $X_r\subset \Pc$ is pfaffian if its 
equation is pfaffian of a skew-symmetric matrix of linear forms in $\Pc$.
\end{definition}

The results proved by Beauville in \cite{beauville} exclude the existence of ACM rank 
$2$ bundles with $c_1=r-1$ on a general hypersurface $X_r$, $r\geq 3$.

\begin{proposition} \label{prop:pfaff} When $r\geq 3$ then a general 
hypersurface $X_r\subset \Pc$ has no normalized indecomposable rank $2$ ACM
bundles $\E$ with $c_1(\E)=r-1$ and $c_2=r(r-1)(2r-1)/6$. 
\end{proposition}
\begin{proof} It follows soon by the following two facts proved in 
\cite{beauville}. $X_r$ is pfaffian if and only if there exists an 
indecomposable ACM rank $2$ vector bundle on $X_r$ with Chern classes
as in the statement. 
Moreover the general hypersurface of degree $r\geq 3$ in 
$\Pc$ is not pfaffian.
\end{proof}

Let us now turn to hypersurfaces of low degree. We want to exclude the existence of 
indecomposable ACM rank $2$ bundles on general hypersurfaces. This follows easily on 
sextic hypersurfaces, using the main result of \cite{cm2}.

\begin{proposition}\label{prop:sextic} On a general sextic hypersurface $X\subset \Pc$ 
all ACM rank $2$ bundles $\E$ split. 
\end{proposition}
\begin{proof} A general hyperplane section $Y$ of $X$ is a general sextic hypersurface 
of $\Pq$. By remark \ref{rm:restrict} we know that an indecomposable ACM rank $2$ bundle 
on $X$ restricts to an indecomposable ACM rank $2$ bundle on $Y$. In \cite{cm2} we 
excluded the existence of such bundles.
\end{proof}

For hypersurfaces of degree $r<6$ we cannot use the same procedure, since there exist 
indecomposable ACM rank $2$ bundles on general cubics, quartics and quintics of $\Pq$.

We use instead an examination of the family of surfaces associated to ACM rank $2$ 
bundles. Let us set some more pieces of notation.

Call $\Hh(d,g)$ the Hilbert scheme of arithmetically Cohen--Macaulay (ACM) surfaces in 
$\Pc$ of degree $d=c_2$ and sectional genus $g$ such that $2g-2=c_2(c_1+r-5)$. This is 
a smooth quasi-projective subvariety of the Hilbert scheme.\par\noindent
Let $\Pj(r)$ be the scheme which parametrizes hypersurfaces of degree  $r$ in $\Pj^5$.
\par\noindent
In the product $\Hh(d,g) \times \Pj(r)$ one has the incidence variety
\begin{equation} \label{eq:inc}
I(r,d,g)=\{ (S,X):X \ \text{is smooth and} \ S \subset X \} \subset
\Hh(d,g) \times \Pj(r)
\end{equation}
with the two obvious projections $p(r):I(r,d,g) \to \Hh(d,g)$ and 
$q(r):I(r,d,g) \to \Pj(r)$.
The fibers of $q(r)$ are projective spaces of fixed dimension (by Riemann-Roch).
\par\noindent

We will show that $I(r,d,g)$ does not dominate $\Pj(r)$ for all choices of $d,g$ 
corresponding to surfaces associated with an indecomposable ACM rank $2$ bundle on a 
general $X_r$. This is achieved in the next sections by computing the dimension of 
$I(r,d,g)$ and observing that it is smaller than $\dim(\Pj(r))$.

Let us see, for instance, what happens for quadric surfaces.

\begin{remark} \label{rm:quad} \rm Any quadric surface $S$ contained in a general 
hypersurface $X_r$, $r\geq 3$, is reduced since $X_r$ contains no planes. Hence it is 
a surface in $\Pt$, that is $S$ is a complete intersection of type $(1,1,2)$ in $\Pc$. 

Thus we may compute the normal bundle $N_S$ of $S$:
$$ 
h^0(N_S) = h^0(\Oc_S(2)\oplus\Oc_S(1)\oplus\Oc_S(1))= 9+4+4 = 17 
$$
hence $\dim(H(2,0))\leq 17$. 
\end{remark}

\begin{proposition} \label{prop:quad} On a general hypersurface $X_r\subset \Pc$ of 
degree $r\geq 3$ there are no indecomposable normalized ACM rank $2$ bundles $\E$ 
with $c_1(\E)=4-r$. \end{proposition}
\begin{proof} 
First we show that any such bundle $\E$ is associated with a complete intersection 
quadric surface.\par\noindent
Consider a global section $s\in H^0(\E)$ 
and its zero-locus $S$. The exact sequence (\ref{eq:kos}) 
here reads  
$$
0 \to \Oc_{X_r} \to \E \to \Ii_{S/X_r}(4-r) \to 0
$$ 
and implies $h^0(\E(r-4))= h^0(\Oc_{X_r}(r-4))$. $h^1(\E(r-4))$ and $h^2(\E(r-4))$ vanish 
by assumptions. By duality $h^4(\E(r-4))=h^0(\E(r-6))=h^0(\Oc_{X_r}(r-6))$. Hence just as 
in proposition \ref{prop:refin} we use Riemann-Roch to prove that $c_2(\E)=2$. So $S$ has 
degree $2$. Since a general $X_r $ contains no planes, $S$ is reduced and the claim is 
proved.\par\noindent 
Let $\Ii_S$ indicate the ideal sheaf of $S$ in $\Pc$. One computes from the resolution 
of $\Ii_S$:
$$h^0(\Ii_S(r))= 2h^0(\Oc_{\Pc}(r-1))- 2h^0(\Oc_{\Pc}(r-3))+ h^0(\Oc_{\Pc}(r-4))$$
and thus one easily sees that:
$$\dim(I(r,2,0))\leq h^0(\Ii_S(r))-1 + h^0(N_S) < \dim(\Pj(r))$$
for $r>2$, which means that the map $q(r)$ above is not dominant. The conclusion follows.
\end{proof}

With the results above we dispose of the case of cubic hypersurfaces:

\begin{proposition} \label{prop:cubic} On a general hypersurface $X:=X_3\subset \Pc$ of 
degree $ 3$ there are no indecomposable ACM rank $2$ bundles. 
\end{proposition}
\begin{proof} By proposition \ref{prop:refin} we know that any normalized indecomposable 
ACM rank 2 bundle on a smooth cubic hypersurface satisfies  $3> c_1(\E)>0$. So only the 
cases $c_1(\E)=1$ and $c_1(\E)=2$ are left. But on a general cubic hypersurface the case 
$c_1(\E)=2$ is excluded by proposition \ref{prop:pfaff} while the case $c_1(\E)=1$ is 
excluded by proposition \ref{prop:quad}.
\end{proof}

\section{Quartic hypersurfaces}

In this section we fix $r=4$. Our goal is to exclude the existence of indecomposable ACM 
rank $2$ bundles $\E$ on a general quartic fourfold $X:=X_4$. We also assume that $\E$ is 
normalized.

Arguing as in proposition \ref{prop:cubic} we know that for such a bundle $\E$ the only 
possibilities for the first Chern classes are $ c_1(\E)=2$ and $c_1(\E)=1$.

We dispose of these cases using a computation for the normal bundle of the zero--locus of 
a general global section of $\E$.

\begin{remark} \label{rm:Gor} \rm If $\E$ is an ACM rank $2$ bundle on a smooth 
hypersurface $X_r$, then the zero--locus $S$ of a global section of $\E$ has 
codimension at most $3$ in $\Pc$. If it has codimension $3$, then it is an  
``arithmetically Gorenstein'' subscheme of codimension $3$ in the 
projective space $\Pc$. 
Thus its ideal sheaf $\Ii_S$ in $\Pc$ has a self dual free resolution of type 
\begin{equation} \label{eq:res}
0 \to \Oc(-e-6) \to \bigoplus_{j=1}^r \Oc(-m_j) \to \bigoplus_{i=1}^r \Oc(-n_i) \to 
\Ii_S \to 0 
\end{equation}
where $e$ is the number such that the canonical class of $S$ is $e$ times the hyperplane 
section and $e+6-m_i=n_i$ for all $i$.
\end{remark}

Using the previous resolution one can compute the cohomology of the normal bundle $N_S$ 
of $S$ in $\Pc$. Indeed by \cite{eb} and Theorem 2.6 of \cite{kl-roi} we have the 
following:

\begin{proposition} \label{KMR} {\rm (Kleppe - Mir\'o-Roig)}
With the notation of the previous remark, order the integers $n_i$ and $m_j$ so that:
$$n_1 \leq n_2 \leq ... \leq n_r \quad\text{ and }\quad m_1 \geq m_2 \geq ... \geq m_r.$$ 
Then:
\begin{equation} \label{eq:kl-roi}
\begin{aligned}
h^0N_S & =\sum_{i=1}^r h^0\Oc_S(n_{i})+
\sum_{1 \leq i<j \leq r} \binom{-n_{i}+m_{j}+5}{5}+\\
& -\sum_{1 \leq i<j \leq r} \binom{n_{i}-m_{j}+5}{5}-
\sum_{i=1}^r \binom{n_{i}+5}{5}.
\end{aligned}
\end{equation}
\end{proposition}

\begin{remark} \label{rm:lift} \rm If $S$ is an ACM subscheme of $\Pc$ and $C$ is a 
general hyperplane section of $S$, then a minimal resolution of the ideal sheaf of $C$ 
in $\Pq$ lifts to a minimal resolution of the ideal sheaf of $S$ in $\Pc$.
\end{remark}

Let us go back to general quartic hypersurfaces $X$.\par
A general hyperplane section $Y$ of $X$ is a general quartic threefold in $\Pq$ and the 
restriction $\E'$ of $\E$ to $Y$ is an indecomposable ACM bundle of rank $2$. These 
bundles are classified in \cite{m2}, where the possibilities for the second Chern classes 
of $\E'$, hence also of $\E$, are listed. These possibilities are:\par
\centerline{$(c_1,c_2)\in\{(-1,1), (0,2), (1,3),(1,4),(1,5),(2,8),(3,14)\}$. }\par
The cases $(c_1,c_2)=(-1,1), (0,2), (3,14)$ cannot occur on a general quartic fourfold, 
by propositions \ref{prop:pfaff}, \ref{prop:refin} and \ref{prop:quad}.\par
 We explore the remaining possibilities  case by case. 

\begin{case} $c_1(\E)=1$, $c_2(\E)=3$.\par\noindent
\rm  By \cite{m2} $\E'$ is associated with a plane cubic curve, hence $\E$ is associated 
with a cubic surface $S\subset\Pt$. It turns out $h^0(N_S) = h^0(\Oc_S(3)\oplus
\Oc^2_S(1)) = 27$ while the ideal sheaf $\Ii_S$ has $h^0(\Ii_S(4))=95$. Thus in this case 
$\dim(I(4,3,1)) \leq 121$. Since $\dim(\Pj(4))=125$, the projection $q(4): I(4,3,1) \to 
\Pj(4)$ cannot be dominant and this case is excluded on a general quartic hypersurface $X$.
\end{case}

\begin{case} $c_1(\E)=1$, $c_2(\E)=4$.\par\noindent
\rm  By \cite{m2} $\E'$ is associated with a quartic curve, complete intersection of $2$ 
quadrics in $\Pt$. Hence $\E$ is associated with a complete intersection of two quadrics 
$S\subset\Pq$. It turns out $h^0(N_S) = h^0(\Oc^2_S(2)\oplus\Oc_S(1)) = 31$ while 
$h^0(\Ii_S(4))=85$, so that $\dim(I(4,4,1)) \leq 115$ and the projection 
$q(4): I(4,4,1) \to \Pj(4)$ cannot be dominant.
\end{case}

\begin{case} $c_1(\E)=1$, $c_2(\E)=5$.\par\noindent
\rm  By \cite{m2} $\E'$ is associated to an elliptic quintic curve, whose ideal sheaf 
$\Ii$ in $\Pq$ has resolution:
\begin{equation}
0 \to \Oc_{\Pq}(-5) \to \Oc^5_{\Pq}(-3) \to \Oc^5_{\Pq}(-2) \to \Ii \to 0
\end{equation}
from which we have the resolution for the ideal sheaf $\Ii_S$ of a quintic surface 
$S$ associated with $\E$. Now we use proposition \ref{KMR} to compute $h^0(N_S)=35$ 
while from the resolution we get $h^0(\Ii_S(4)) = 75$ so that $\dim(I(4,5,1))\leq 109$ 
and again $q(4)$ does not dominate $\Pj(4)$.
\end{case}


\begin{case} $c_1(\E)=2$, $c_2(\E)=8$.\par\noindent
\rm  By \cite{m2} $\E'$ is associated to a curve of degree $8$ in $\Pq$, whose ideal 
sheaf $\Ii$ has resolution:
\begin{equation}
0 \to \Oc_{\Pq}(-6) \to \Oc^3_{\Pq}(-4)\oplus\Oc^x_{\Pq}(-3)\to \Oc^x_{\Pq}(-3)\oplus
\Oc^3_{\Pq}(-2) \to \Ii \to 0
\end{equation}
from which we have the resolution for the ideal sheaf $\Ii_S$ of a surface $S$ of degree 
$8$ associated with $\E$. Notice that we do not know the number of minimal generators of 
degree $3$ for the ideal sheaf of $S$ (if any). Nevertheless we may use proposition 
\ref{KMR} to compute $h^0(N_S)$. Indeed in the computation it turns out that the 
contribution of cubic generators disappears and one gets $h^0(N_S)=54$. Also one sees 
that $h^0(\Ii_S(4)) = 60$ so 
that $\dim(I(4,8,5))\leq 113$ and again $q(4)$ does not dominate $\Pj(4)$.
\end{case}

No other cases may occur, by \cite{m2}. Hence we conclude:

\begin{proposition} \label{prop:quartic} On a general hypersurface $X\subset \Pc$ of 
degree $4$ there are no indecomposable ACM rank $2$ bundles. 
\end{proposition}

\section{Quintic hypersurfaces}

In this section we exclude the existence of indecomposable ACM rank $2$ bundles $\E$ on 
a general quintic fourfold $X$. As usual we assume that $\E$ is normalized.

In this case, we are left with several cases for the first Chern class, namely $c_1(\E)=
0,1,2,3$.

Again a general hyperplane section $Y$ of $X$ is a general quintic threefold in $\Pq$ 
and the restriction $\E'$ of $\E$ to $Y$ is an indecomposable ACM bundle of rank $2$. 
These bundles are classified in \cite{cm1}. In particular for the Chern classes we have 
the following possibilities:

$$
\begin{matrix}
c_1 & c_2 \\
0 & 3,4,5 \\
1 & 4,6,8 \\
2 & 11,12,13,14 \\
3 & 20
\end{matrix}
$$

We explore again the situation case by case. 

\begin{case} $c_1(\E)=0$, $c_2(\E)=3$.\par\noindent
\rm  By \cite{cm1} $\E'$ is associated to a plane cubic curve, hence $\E$ is associated 
to a cubic surface $S\subset\Pt$. Then as above one computes $h^0(N_S) = 27$ while
the ideal sheaf $\Ii_S$ has $h^0(\Ii_S(5))=206$. Thus in this case 
$\dim(I(5,3,1)) \leq 232$. Since $\dim(\Pj(5))=251$, the projection $q(5): I(5,3,1) \to 
\Pj(5)$ cannot be dominant and this case is excluded on a general quintic hypersurface.
\end{case}

\begin{case} $c_1(\E)=0$, $c_2(\E)=4$.\par\noindent
\rm  By \cite{cm1} $\E'$ is associated with a quartic curve, complete intersection of 
$2$ quadrics in $\Pt$. Hence $\E$ is associated to a complete intersection of two quadrics 
$S\subset\Pq$. It turns out $h^0(N_S) = 31$ while $h^0(\Ii_S(5))=191$, so that  
$\dim(I(5,4,1)) \leq 221$ and $q(5)$ is not dominant.
\end{case}

\begin{case} $c_1(\E)=0$, $c_2(\E)=5$.\par\noindent
\rm  By \cite{cm1} $\E'$ is associated with a quintic elliptic curve and as above 
one gets a resolution
\begin{equation}
0 \to \Oc_{\Pc}(-5) \to \Oc^5_{\Pc}(-3) \to \Oc^5_{\Pc}(-2) \to \Ii \to 0
\end{equation}
for the ideal sheaf of a surface associated with $\E$. Then $h^0(N_S) = 35$ while 
$h^0(\Ii_S(5))=176$, so that  $\dim(I(5,5,1)) \leq 210$ and $q(5)$ is not dominant.
\end{case}

\begin{case} $c_1(\E)=1$, $c_2(\E)=4$.\par\noindent
\rm  By \cite{cm1} $\E'$ is associated with a plane quartic curve, hence $\E$ is 
associated with a quartic surface $S\subset\Pt$. It turns out $h^0(N_S) = h^0(\Oc_S(4)
\oplus\Oc^2_S(1)) = 42$ while the ideal sheaf $\Ii_S$ has $h^0(\Ii_S(5))=200$. Thus in 
this case $\dim(I(5,4,3)) \leq 241$ and the projection $q(5): I(5,4,3) \to \Pj(5)$ cannot 
be dominant.
\end{case}

\begin{case} $c_1(\E)=1$, $c_2(\E)=6$.\par\noindent
\rm  By \cite{cm1} $\E'$ is associated with a sextic curve, complete intersection of a 
quadric and a cubic in $\Pt$. Hence $\E$ is associated with a corresponding complete 
intersection  $S\subset\Pq$. It turns out $h^0(N_S) = 48$ while $h^0(\Ii_S(5))=175$, 
so that  $\dim(I(5,6,4)) \leq 222$ and the projection $q(5)$ cannot be dominant.
\end{case}

\begin{case} $c_1(\E)=1$, $c_2(\E)=8$.\par\noindent
\rm  By \cite{cm1} $\E'$ is associated to a curve of degree $8$ in $\Pq$, whose ideal 
sheaf $\Ii$ has resolution:
\begin{equation}
0 \to \Oc_{\Pq}(-6) \to \Oc^3_{\Pq}(-4)\oplus\Oc^x_{\Pq}(-3)\to \Oc^x_{\Pq}(-3)\oplus
\Oc^3_{\Pq}(-2) \to \Ii \to 0.
\end{equation}
As above one gets $h^0(N_S)=54$ while $h^0(\Ii_S(5)) = 150$ so that 
$\dim(I(5,8,5))\leq 203$ and  $q(5)$ is not dominant.
\end{case}

Consider now the case $c_1(\E)=2$ and $c_2(\E)= 11,12,13,14$. Let $S$ be a surface 
associated with $\E$ and call $C$ a general hyperplane section of $S$, which is thus 
associated with $\E'$. One computes:
$$  h^0(\Ii_S(5))= 245 - 10 \deg(S) =245 -10 c_2(\E)$$
so that we only need to prove that:
\begin{equation}\label{Ns}
h^0(N_S)<10c_2(\E)+7.
\end{equation}  

We use the results of \cite{cm1} \S 4 and \cite{cm2} case 5.7 to compute a minimal 
resolution for the ideal sheaf of $C$ in $\Pq$, hence also a resolution of $\Ii_S$, 
which leads to the computation of $h^0(N_S)$, via proposition \ref{KMR}.

\begin{case} $c_1(\E)=2$, $c_2(\E)=11$.\par\noindent
\rm  By \cite{cm1} \S 4 the resolution of the ideal sheaf $\Ii_S$ is:
\begin{equation}
\begin{aligned}
0 & \to \Oc_{\Pc} (-7) \to \Oc^b_{\Pc}(-3) \oplus \Oc^c_{\Pc}(-4) \oplus \Oc^3_{\Pc}(-5) 
\to \\
& \to \Oc_{\Pc}^3(-2)\oplus \Oc_{\Pc}^c(-3) \oplus \Oc_{\Pc}^b(-4) \to \Ii_S \to 0.
\end{aligned}
\end{equation}
Comparing the first Chen classes in the exact sequence, one finds $c=b-2$. Using equation 
\ref{eq:kl-roi} one is able to compute $h^0(N_S)$. It turns out that $b$ and $c$ cancel 
and one finds $h^0(N_S)=83<117$ so that $\dim(I(5,11,12))\leq 214$ and  $q(5)$ is not 
dominant.
\end{case}

\begin{case} $c_1(\E)=2$, $c_2(\E)=12$.\par\noindent
\rm  By \cite{cm1} \S 4 the resolution of the ideal sheaf $\Ii_S$ is:
\begin{equation}
\begin{aligned}
0 & \to \Oc_{\Pc}(-7) \to \Oc_{\Pc}^b(-3) \oplus \Oc_{\Pc}^c(-4) \oplus \Oc_{\Pc}^2(-5) 
\to \\
& \to \Oc_{\Pc}^2(-2) \oplus \Oc_{\Pc}^c(-3) \oplus \Oc_{\Pc}^b(-4) \to \Ii_S \to 0
\end{aligned}
\end{equation}
where $b=c-1$ and $b=0,1$, according with the existence of a cubic syzygy between the 
two quadrics. In both cases, using equation \ref{eq:kl-roi} one computes $h^0(N_S)=81<127$ 
so that $\dim(I(5,12,13))\leq 205$ and  $q(5)$ is not dominant.
\end{case}

\begin{case} $c_1(\E)=2$, $c_2(\E)=13$.\par\noindent
\rm  In this case we have only one quadric containing $S$ and the resolution of $\Ii_S$ 
is given by:
\begin{equation}
\begin{aligned}
0 & \to \Oc_{\Pc} (-7) \to \oplus \Oc_{\Pc}^4(-4)\oplus \Oc_{\Pc}(-5) \to \\
& \to \Oc_{\Pc}(-2) \oplus \Oc_{\Pc}^4(-3) \to \Ii_S \to 0.
\end{aligned}
\end{equation}
So one computes $h^0(N_S)=79<137$. It follows that $q(5)$ is not dominant.
\end{case}

\begin{case} $c_1(\E)=2$, $c_2(\E)=14$.\par\noindent
\rm  By \cite{cm1} \S 4 the resolution of the ideal sheaf $\Ii_S$ is:
\begin{equation}
0 \to \Oc_{\Pc} (-7) \to \oplus \Oc_{\Pc}^7(-4) \to \oplus \Oc_{\Pc}^7(-3) \to \Ii_S \to 0
\end{equation}
and one computes $h^0(N_S)=77<147$ so that $q(5)$ is not dominant.
\end{case}

Finally for $c_1=3$ we have:

\begin{case} $c_1(\E)=3$, $c_2(\E)=20$.\par\noindent
\rm  By \cite{cm1} we know the resolution of the ideal sheaf of a curve associated with 
$\E'$, so that the ideal sheaf of a surface associated with $\E$ is:
\begin{equation}
0 \to \Oc_{\Pc}(-8) \to \Oc^4_{\Pc}(-5)\to \Oc^4_{\Pc}(-3) \to \Ii_S \to 0.
\end{equation}
One computes $h^0(N_S)=110$ and $h^0(\Ii_S(5)) = 80$ so that $\dim(I(5,20,31))\leq 189$ 
and  $q(5)$ is not dominant.
\end{case}

Hence we may conclude

\begin{proposition} \label{prop:quintic} On a general hypersurface $X\subset \Pc$ of 
degree $5$ there are no indecomposable ACM rank $2$ bundles. 
\end{proposition}

The main theorem follows.

\begin{remark}\rm
By \cite{bgs} there exists a non discrete family (up to twist) of isomorphism classes of 
indecomposable ACM vector bundles on any smooth projective hypersurface of degree $r 
\geq 3$ in the $5$--dimensional complex projective space $\Pj^5$. On a general $X_r$ the 
rank of the bundles constructed in \cite{bgs} is $16$ (cfr. \cite{m3}). \par
The problem of determining the minimum rank $BGS(X_r)$ for ACM bundles on $X_r$ moving in 
a non--trivial family (the {\it BGS invariant)} is still open.\par
We prove in this paper that $BGS(X_r)>2$ for general hypersurfaces in $\Pc$ of degree 
$r\leq 6$.  
\end{remark}


\begin{thebibliography}{ACGH}

\bibitem{arrondo} E.Arrondo,
\emph{On congruences of lines in the projective space},
M\'em. Soc. Math. France (N.S.) No. \textbf{50}\ (1992), 96 pp.

\bibitem{AC} E.Arrondo and L.Costa,
\emph{Vector bundles on Fano 3-folds without intermediate cohomology},
Comm. Algebra \textbf{28} (2000), no. 8, 3899--3911.

\bibitem{beauville} A.Beauville,
\emph{Determinantal hypersurfaces. 
Dedicated to William  Fulton on the occasion of his 60th birthday},
Michigan Math. J.  \textbf{48} (2000), 39--64.

\bibitem{bgs} R.O.Buchweitz, G.M.Greuel, and F.O.Schreyer,
\emph{Cohen-Macaulay modules on hypersurface singularities II},
Invent. Math. \textbf{88} (1987), 165--182.

\bibitem{cm1} L.Chiantini and C.Madonna,
{\it ACM bundles on a general quintic threefold. 
Dedicated to Silvio Greco on the occasion of his 60th birthday (Catania, 2001)},  
Matematiche (Catania) \textbf{55} (2000), no. 2, 239--258.

\bibitem{cm2} L.Chiantini and C.Madonna,
{\it A splitting criterion for rank $2$ bundles on a general sextic threefold},
Internat. J. Math. \textbf{15} (2004), no. 4, 341--359.

\bibitem{debarre} O.Debarre, 
\emph{Sur la vari\'et\'e des espaces lin\'eaires contenus dans une intersection 
compl\'ete}, Math. Ann. \textbf{312} (1998), no. 3, 549--574.

\bibitem{EG} E. G. Evans and P. Griffith, {\it The syzygy problem}, Ann. of Math. 
\textbf{214}  (1981), 323--333.  

\bibitem{eb} D. Buchsbaum and D. Eisenbud, 
{\it Algebra structures for finite free resolutions, and some structure theorems for 
ideals of codimension 3},  Amer. J. Math. \textbf{99} (1977) no.3 447--485.

\bibitem{H1} R.Hartshorne,
{\it Stable vector bundles of rank 2 on $\mathbb{P}^3$},
Math. Ann. \textbf{238} (1978), 229--280.

\bibitem{kl-roi}  J.O.Kleppe and R.M.Mir\'o--Roig,
{\it The dimension of the Hilbert scheme of Gorenstein codimension 3 subschemes}, 
J.Pure Appl. Alg. \textbf{127} (1998),73--82.

\bibitem{m0} C.Madonna, 
{\it A splitting criterion for rank 2 vector bundles on hypersurfaces in $ P\sp 4$}, 
Rend. Sem. Mat. Univ. Politec. Torino  \textbf{56} (1998),  no. 2, 43--54.

\bibitem{m2} C.Madonna, 
{\it Rank-two vector bundles on general quartic hypersurfaces in ${\mathbf P}\sp 4$},
Rev. Mat. Complut. \textbf{13} (2000),  no. 2, 287--301.

\bibitem{m3} C.Madonna, {\it Rank 4 ACM bundles on a  smooth quintic hypersurface in 
$\Pj^4$}, preprint.

\bibitem{OSS} C.Okonek, M.Schneider and H.Spindler,
Vector bundles on complex projective spaces, 
Progress in Mathematics \textbf{3}, 1980.

\bibitem{Ott} G.Ottaviani, 
{\it Some extension of Horrocks criterion to vector bundles on Grassmannians and 
quadrics}, Ann. Mat. Pura Appl. (4) \textbf{155} (1989), 317--341.

\bibitem{v} M.Roggero and P.Valabrega, 
{\it The  speciality lemma, rank 2 bundles and Gherardelli-type theorems for 
surfaces in $\Pj^4$}, Compositio Math. \textbf{139} (2003),  no. 1, 101--111.

\end{thebibliography}
\end{document}